\documentclass[11pt]{article}
\usepackage{graphicx}
\usepackage{authblk}
\usepackage{mathtools}
\usepackage{amsmath}
\usepackage{amsthm}
\newcommand\lar[2]
  {\expandafter\gdef\csname labeled:#1\endcsname{#2}%
   \label{#1}#2}
\newcommand\recallLabel[1]
   {\csname labeled:#1\endcsname\tag{\ref{#1}}}
   
\usepackage{amssymb, amsmath}
\usepackage{epstopdf, wrapfig, framed,caption}
\usepackage{color}
\usepackage[makeroom]{cancel}
\usepackage{stmaryrd}
\usepackage{hyperref}
\usepackage{makecell}
\usepackage{authblk}
\usepackage[extra]{tipa}
\newcommand{\Frac}[2]{\frac{\textstyle #1}
                           {\textstyle #2}}

\def\mbbI{\mathbb I}

\def\op{\overline{p}}
\def\oq{\overline{q}}
\def\ou{\overline{U}}
\def\ov{\overline{W}}
\def\oomega{\overline{\omega}}

\newcommand{\eps} {\epsilon}

\newcommand{\be}{\begin{equation}}
\newcommand{\ee}{\end{equation}}
\newcommand{\ba}[1] {\begin{array}{ #1 }}
\newcommand{\ea}{\end{array}}
\def\tu{ \tilde{U} }
\def\tU{\tilde{U}}

\def\mrL{\mathrm L}
\def\tE{\tilde{\mathrm E}}
\def\oE{\overline{\mathrm E}}
\def\tx{\tilde{x}}
\def\tF{\tilde{F}}
\def\tR{\tilde{R}}

\def\mrD{\mathrm D}
\def\mrC{\mathrm C}

\def\tp{\tilde{p}}
\def\tq{\tilde{q}}

\newcommand{\sech} { {\rm sech} \hskip 0.01in}

\newtheorem{thm}{Theorem}
\newtheorem{lemma}[thm]{Lemma}

\newtheorem{prop}[thm]{Proposition}

\renewcommand{\appendix}{\Alph{section}}
\title{Curvature and Chaos  in the Defocusing Parameteric Nonlinear Schr\"odinger System}

\author[1]{Keith Promislow}
\author[2]{Abba Ramadan}
\affil[1]{Department of Mathematics, Michigan State University,
East Lansing, MI 48824, USA}
\affil[2]{Department of Mathematics,
The University of Alabama, 
Tuscaloosa, AL 35401, USA}

\begin{document}
\maketitle
\begin{abstract}
The parameteric nonlinear Schr\"odinger equation models a variety of parametrically forced and damped dispersive waves. For the defocusing regime, we derive a normal velocity for the evolution of curved dark-soliton fronts that represent a $\pi$-phase shift across a thin interface.  We establish that depending upon the strength of parametric term the normal velocity evolution can transition from a curvature driven flow to motion against curvature regularized by surface diffusion of curvature. In the former case interfacial length shrinks, while in the later the interface length generically grows until self-intersection followed by a transition to chaotic motion.
\end{abstract}
 


\section{Introduction}

The parametric nonlinear Schr\"odinger (PNLS) equation is a general model for parametrically forced surface waves and for pattern formation.  It has been derived in the context of Faraday waves \cite{bib:PRL95} where increased driving force drives transitions to zigzag patterns and chaotic behavior. The  PNLS has also be derived as a model of phase sensitive amplifiers \cite{bib:AGJS} and in the large detuning limit of optical parametric oscillator systems \cite{bib:PK, bib:OPO1, bib:OPO2}.  More recently it has been proposed as a model for dissipative self organization \cite{bib:OPO4} and as a template for second-order phase transitions between degenerate and non-degenerate regimes, \cite{bib:OPO3}.   

We present an analysis of the evolution of curved dark soliton fronts in the 1+2D PNLS equation. These fronts represent $\pi$-phase shifts in the optical field. We  consider interfaces that have bounded curvatures and derive a normal velocity that describes the temporal evolution of the front. In particular if the parametric strength decreases through a critical value the normal velocity transitions from a curvature driven flow to motion \emph{against} curvature regularized by surface diffusion of curvature. Specifically in the limit in which the ratio of dispersive length scale to domain size is small ($\eps\ll 1$) we identify a bifurcation parameter $\mu$ for which the normal velocity of the dark soliton interface admits the expansion
\be \label{e:dsV}
V= -\alpha_0 \kappa_0 +\eps^2(\nu\Delta_s \kappa_0 +\zeta \kappa_0^3) +O(\eps^4).
\ee 
Here $\kappa_0$ is the curvature of the interface, $\Delta_s$ is the Laplace-Beltrami surface diffusion operator, and $\alpha_0$, $\nu$ and $\zeta$ are $\mu$-dependent real coefficients. The leading order coefficient $\alpha_0$ has the same sign as $\mu$, its sign change encodes the transition between motion with and against curvature. Crucially we establish the existence of $\mu_*>0$, independent of $\eps$, such that $\nu>0$ if $|\mu|<\mu_*.$ This allows the surface diffusion to regularize the motion against curvature that arises for $\mu<0.$ Such curvature flow transitions have been studied in dissipative systems such as polymer melts, \cite{bib:CP22} but their presence in a dispersive system are here-to-for unstudied. The analysis presented here is formal but is complemented with a sharp characterization of the transverse spectrum of the wave which shows that the curvature transition is not associated to any transverse instability of the dark soliton.

The paper is organized as follows. In section 2 the PNLS system and the analysis of the 1D transfer spectral problem are presented. In section 3 the inner and outer asymptotic formulations of the system are developed and normal velocity is resolved. Section 4 presents consequences of the normal velocity and compares these to simulations of the full system.

\section{PNLS and One-dimensional Spectral Analysis}
The PNLS system describes the evolution of a complex phase $\Theta\in H^2(\Omega)$,
\be i\Theta_t +\frac{\eps^2}{2}\Delta \Theta-|\Theta|^2\Theta+(i+a)\Theta-\gamma\Theta^*=0,
\label{dfpnls}\ee
on a spatially periodic domain $\Omega=[0,L]^2.$ Here $0<\eps\ll1$ is a small parameter that characterizes the ratio of the dispersive lengthscale to domain size, $a$ is a phase rotation and $\gamma$ is the parametric pump strength.
We rescale the complex phase $\Theta$, time $t$
\be \begin{aligned}
\Theta&=\frac{2}{\sqrt{\beta}} u(\tx)e^{i\theta},\\
\tau&=2 t/\beta,
\end{aligned}
\ee
and space $\tx=2 x/\sqrt{\beta}$ where
$$\beta:=\frac{4}{a+\sqrt{\gamma^2-1}}>0,$$ and the phase angle $\theta$ as the solution of
$$\gamma e^{-2i\theta}=-\sqrt{\gamma^2-1}+i.$$ 
We drop the tilde's and introduce the real vector function $U=(\Re u, \Im u)^t$ which satisfies
\be 
U_\tau=F(U) := \begin{pmatrix}
0 & -(\eps^2\Delta -2|U|^2+1-\mu)\cr 
(\eps^2\Delta -2|U|^2+2) & -\beta
\end{pmatrix} U,\label{e:FMU}\ee
where the bifurcation parameter
$$ \mu :=-\Frac{a-3\sqrt{\gamma^2-1}}{a+\sqrt{\gamma^2-1}},$$
lies in $[-1,3]$. In what follows, a shift from $\mu>0$ to $\mu<0$ will trigger the transition in curvature motion.

Posed on the line ${\mathbb R}$, in the scaled coordinate $z=x/\eps$, and the PNLS system has a dark-soliton steady-state solution
\be \Phi_0(z)=\begin{pmatrix}
    \phi(z) \cr 0 
\end{pmatrix},\label{front}\ee
where $\phi(z)=\tanh(z),$ solves
\be \label{e:front-eq} \partial_z^2\phi -2\phi^3+2\phi=0.
\ee
The linearization of the 1D PNLS system about $\Phi_0$ yields the system
\be W_\tau = \mrL W,\label{lineq}\ee
where the 1D linear operator $\mrL$ is given by
\be \mrL = \begin{pmatrix} 0 & \mrD \cr -\mrC& -\beta\end{pmatrix},\label{linop}\ee
with
\be \ba{rcl}
  \mrC &=& -\partial_z^2 -6 \psi^2 +4,\\
  \mrD &=& -\partial_z^2 -2 \psi^2  +\mu +1,\ea\ee
and $\psi(x)=\sech(x).$
These operators have point spectrum eigenfunction pairs
\be\ba{rcl} \sigma_p(\mrC)&=&\{ (0,\phi^\prime), 
                    (3,\phi\psi)\},\\
\sigma_p(\mrD)&=&\{(\mu,\psi)\}.\ea\label{CDspec}\ee
The operators $\mrC$ and $\mrD$ may have other point spectrum in their respective gaps $(3,4)$ and $(\mu,1+\mu)$ between their largest point spectrum and the branch point of their essential spectrum. The characterization of these ground state eigenvalues and the essential spectrum show that $\mrC\geq0$ and $
\mrD>0$ if $\mu>0$, while the dimension  $n(\mrD)$ of the negative space of $\mrD$ satisfies $n(\mrD)\geq 1$ if $\mu<0.$ 
A structural point of the analysis arises from the generic fact that that ground states of the Sturmian operators $C$ and $D$ are both non-zero with full support, and hence can not be orthogonal. 
 
\subsection{Essential Spectrum of $\mrL$}
To characterize the essential spectrum of $\mrL$ substitute $W=e^{\lambda\tau}e^{isx}V$ into \eqref{lineq} to obtain an eigenvalue problem for $V,$
\be \lambda V=\begin{pmatrix} 0 & s^2+2(1+\mu)\cr
             -s^2-4 & -\beta\end{pmatrix}V,\ee
which has nontrivial solutions for
\be \lambda=\Frac{-\beta\pm \sqrt{\beta^2-4(s^2+1+\mu)
(s^2+4)}}{2}.\ee
The maximum of $\Re\lambda$ occurs at $s=0$, and hence all
$\lambda\in\sigma_{\rm ess}(\mrL)$ satisfy
\be \Re\lambda \leq \lambda_M :=\Re\left(\Frac{-\beta+
\sqrt{\beta^2-16(1+\mu)}}{2}\right)<0.\ee

\subsection{The Kernel of $\mrL$}
 The operator $C$ has a kernel and from Lemma 3.5 of \cite{bib:PK}  we know that $\lambda=0$ is a simple eigenvalue of $\mrL$, for all $\mu.$ 
For $\mu\neq 0$ the kernel of $\mrL$ and $\mrL^\dag$ are spanned by
\be \Psi_0 = \begin{pmatrix}\phi^\prime \cr 0\end{pmatrix} \quad\quad
    \Psi_0^\dag = \begin{pmatrix}\beta \mrD^{-1}\phi^\prime \cr \phi^\prime\end{pmatrix}.
\label{psi0}\ee
 For $\mu=0$ the kernel remains simple and
 \be \Psi_0 = \begin{pmatrix}\phi^\prime \cr 0\end{pmatrix} \quad\quad
    \Psi_0^\dag = \begin{pmatrix}\psi \cr 0 \end{pmatrix}.
\label{psi0mu0}\ee
When written with unit norm, the eigenfunctions are smooth in $\mu.$

The inverse of $\mrL$ is given by
\be \label{e:Linv-gen} \mrL^{-1}=\begin{pmatrix}-\beta \mrC^{-1} \mrD^{-1} & -\mrC^{-1}\cr
                        \mrD^{-1} & 0\end{pmatrix}.\ee
When $C$ and $D$ have kernels care is required to insure that the inverses act on their domain.

\subsection{Point Spectrum of $\mrL$}

Theorem 3.6 of
\cite{bib:PK} establishes the existence of $\lambda_M>0$ such that
$$ \sigma(\mrL)\subset \{\lambda \bigl| \Re\lambda<-\lambda_M\}\cup\{0\}.$$ We provide an alternate proof that readily generalizes to accommodate in-plane modes.

The point spectrum of $\mrL$ 
is comprised of eigenfunctions localized in $x$ that solve
\be \lambda P = \begin{pmatrix} 0 &\mrD\cr -\mrC & -\beta\end{pmatrix}P.\ee
As two equations for the two unknowns $P=(P_1,P_2)^t$ satisfies
\be\ba{rcl} \lambda P_1& =&DP_2,\\
            (\lambda+\beta) P_2&=& -\mrC P_1.\ea\label{ptspectrumeq}\ee
Assuming that $\lambda\neq 0$ then $P\perp \Psi_0^\dag.$  For $\mu\neq 0$ we combine this with the first equation of \eqref{psi0} and deduce that either $\lambda=-\beta<0$ or
$$ P_1\perp \mrD^{-1}\phi' \hspace{0.5in}\textrm{and}\hspace{0.5in} P_2\perp \phi'.$$
For $\mu\neq 0$ the operator $D>0$ is invertible and  we may solve for $P_1$,
$$\mrC P_1 +\lambda(\lambda+\beta)\mrD^{-1} P_1=0.$$
Taking the complex-valued inner product with $P_1$ yields the relation
\be 
\label{e:def_rho1}
\lambda(\lambda+\beta) = -\frac{\langle CP_1,P_1\rangle}{\langle \mrD^{-1}P_1,P_1\rangle}=:\rho_1.
\ee
The quadratic formula shows that
\be 
\label{e:spec_cross}
\lambda= \frac{-\beta \pm \sqrt{\beta^2-4\rho_1}}{2}.
\ee
In particular we deduce that $\sigma_p(\mrL)$ resides in the range of the right-hand side over the possible values of $\rho_1.$  
In particular if $\rho_1>0$ then $\Re\lambda<-\rho_1/\beta<0.$ 
For $\mu\neq 0$ this motivates the definition of $X_*(\mu)=\{\mrD^{-1}\phi'\}^\bot$ and the real number,
\be
\label{e:def_rho*}
\rho_*(\mu):=\min_{P_1\in X_*} \frac{\langle \mrC P_1,P_1\rangle}{\langle \mrD^{-1}P_1,P_1\rangle}. 
\ee

\begin{lemma}
  There exists $\mu_*,d_+>0$ such that $D\big|_{X_*}>d_+$ for all $\mu\in[-\mu_*,3].$
\end{lemma}
\begin{proof}
We apply Proposition 5.3.1 of \cite{bib:KP13} to the operator $D$ constrained to act on $X_*$. Taking $\mu_*<0$ with $|\mu_*|$ sufficiently small, then $D$ has negative index $n(\mrD)\leq1$ for all $\mu\in[-\mu_*,3].$ We deduce that
\be \label{e:D-index} n\left(D\bigl|_{X_*}\right)=n(\mrD)-n(A),\ee
where $A=\langle \mrD^{-1}(\mrD^{-1}\phi'),\mrD^{-1}\phi'\rangle\in{\mathbb R}.$ Recalling that $\psi$ is the ground state of $D$ with eigenvalue $\mu$, we write 
\be 
\label{e:phi-bot}
\phi'=\frac{\langle \phi',\psi\rangle}{\lVert\phi'\rVert\lVert\psi\rVert} \psi +\psi^\bot,
\ee
where $\psi^\bot\in X_D:=\{\psi\}^\bot$ satisfies $\lVert\psi^\bot\rVert\leq \lVert\phi'\lVert.$  In particular we have the relation
$$D^{-3} \phi' = \frac{\langle \phi',\psi\rangle}{\Vert\phi'\Vert\Vert\psi\Vert}\frac{ \psi}{\mu^3}+D^{-3}\psi^\bot.$$
Since $X_D$ is a spectral subspace of $D$ it follows that $\sigma\left(D\big|_{X_D}\right)=\sigma(\mrD)\backslash\{\mu\}.$ Hence there exists $\mu_*>0$ and a constant $\tilde d_+>0$ such that $D\big|_{X_D}\geq \tilde d_+$ for $\mu\in[-\mu_*,3].$ We deduce that 
$$ \left |A - \frac{|\langle \phi',\psi\rangle|^2}{\Vert\phi'\Vert\Vert\psi\Vert} \frac{1}{\mu^3}\right|\leq \tilde d_+^{-3}\Vert\phi'\Vert^2.$$
It follows that $A<0$ if $\mu\in[-\mu_*,0)$  for $\mu_*>0$ sufficiently small. Moreover we have the limit $A\to-\infty$ at $\mu\to0^-.$ The index relation \eqref{e:D-index} implies that $n\left(D\big|_{X_*}\right)=0$ for $\mu\in(-\mu_*,0].$ Since $D>0$ for $\mu>0$,
the negative index of $D$ is zero for $\mu\in[-\mu_*,3].$ The lower bound of $D\big|_{X_*}$ is given by its ground state eigenvalue. The ground state eigenvalue is continuous in $\mu,$ and the existence of $d_+>0$ follows. 
\end{proof}

\begin{prop}
 There exists $\mu_*>0$ such that $\rho_*(\mu)>0$ for $\mu\in[-\mu_*,3]$.   
\end{prop}

\begin{proof}
    The operator $C$ has a simple kernel spanned by $\phi'$ and is strictly positive on $\{\phi'\}^\perp.$ The operator $C\big|_{X_*}$ is strictly positive so long as
    $$\langle \mrD^{-1}\phi', \phi'\rangle\neq 0.$$ 
 This is obvious for $\mu>0$ since $D$ is positive there. For $\mu<0$ we use the decomposition \eqref{e:phi-bot} to write
\be 
\label{e:A1}
\left| \langle \mrD^{-1}\phi', \phi'\rangle-
    \frac{\langle \phi',\psi\rangle^2}{\Vert\phi'\Vert\Vert\psi\Vert} \frac{1}{\mu}\right|\leq \tilde d_+\Vert\phi'\Vert^2.
\ee
  This implies the existence of $\mu_*>0$ for which the inner product is not zero for all $\mu\in[-\mu_*,3].$
  That is there exists a constant $c_+>0$ such that $\mrC\big|_{X_*}\geq c_+$ for all values of $\mu\in[-\mu_*,3].$ 
  For these $\mu$ we deduce that for all $P_1\in X_*$
  $$
  \begin{aligned}
  \langle \mrC P_1,P_1\rangle &\geq c_+\Vert P_1\Vert^2, \\
  \langle \mrD^{-1} P_1, P_1\rangle &\leq d_+^{-1} \Vert P_1\Vert^2.
  \end{aligned}
  $$
  It follows that $\rho_*> c_+d_+>0$ for these $\mu.$
\end{proof}

\begin{thm}
    There exists $\mu_*,\lambda_M>0$ such that for all $\mu\in[-\mu_*,3]$ the spectrum of $\mrL$ satisfies 
    \be \sigma(\mrL)\subset\{0\}\cup \{\Re\lambda<-\lambda_M\}.
    \ee
Moreover the kernel of $\mrL$ is simple.
\end{thm}
\begin{proof}
The essential spectrum of $\mrL$ lies strictly in the left-half complex plane. If $\lambda\in\sigma_p(\mrL)\backslash\{0\},$ then the relation \eqref{e:spec_cross} holds with $\rho_1$ as defined in \eqref{e:def_rho1}. By definition of $\rho_*$ we have $\rho_1>\rho_*>0$. From the Taylor expansion of the right-hand side of \eqref{e:spec_cross}
$$\Re \lambda<-\frac{\rho_*}{\beta}. $$
Defining $\lambda_M=\rho_*/\beta$ completes the proof.
\end{proof}

\section{Curvature Driven Flow}

We consider a smooth, closed interface 
$\Gamma=\{\gamma(s)\,\big| s\in[0,L]\}$ that breaks $\Omega$ into two regions $\{\Omega_+,\Omega_-\}$ and introduce the local Frenet coordinates
\be\label{e:Frenet}
x= \gamma(s) + n(s) z/\eps,
\ee
where $n(s)$ is the unit outward normal to the curve $\Gamma$ at point $\gamma(s)$ and
 $z$ is signed, $\eps$-scaled distance to $\Gamma.$
If the interface $\gamma$ has smooth, bounded curvatures and is far from self-intersection then the change of variables from $x=(x_1,x_2)$ to $(s,z)$ is well defined on a neighborhood of $\Gamma$. Indeed there exists $\ell>0$ such that the neighborhood contains all points $x\in\Omega$ for which the scaled distance to $\Gamma$ satisfies $z(x)<\ell/\eps.$ We introduce $\overline{\phi}$  which is a smooth function that agrees with $\phi$ for $|z|<\ell/(2\eps)$ and is
identically $1$ for $z>\ell/\eps$ and identically $-1$ for $z<-\ell/\eps.$ The truncated function $\overline{\phi}$ to induces a smooth function $\Phi$ defined on $\Omega$,
\be 
\label{e:Ansatz}
\Phi(x):=\begin{pmatrix} \overline\phi(z(x)) \cr 0 \end{pmatrix}\hspace{0.5in} |z|<\ell/\eps, 
\ee
and $\Phi(x)=\pm1$ if $z>\ell/\eps$ or $z<-\ell/\eps$ respectively. Since $\phi$ decays exponentially to constant values at an $O(1)$ rate in $z$, this modification induces exponentially small perturbations the do not impact the analysis. The overbar on $\phi$ is dropped in the sequel.  The evolution of $U$ is tracked via its interface map $\gamma=\gamma(s,\tau)$ whose motion prescribed by the normal velocity which controls the evolution of the curvature through the relation
\eqref{e:Pismen}.
Knowledge of the curvatures is equivalent to prescribing $\gamma$ up to rigid body motion.

\subsection{Outer Expansion}
The  outer region is divided into inside $z<0$ and outside $z>0$ sets, $\Omega_\pm.$ These regions are described by Cartesian variables. We expand the ansatz as 
$$U=\begin{pmatrix}p \cr q \end{pmatrix}= u_0+\eps u_1+\eps^2 u_2+O(\eps^3),
$$ 
where each term has a vector decomposition 
$$u_i=\begin{pmatrix}p_i \cr q_i \end{pmatrix}, \hspace{0.5in}
i= 0, 1, 2, \ldots.$$
To match with the ansatz \eqref{e:Ansatz} we impose
$$u_0=\begin{pmatrix} p_0\cr q_0\end{pmatrix}=\begin{pmatrix}  \mbbI_{\Omega_+}-\mbbI_{\Omega_-} \cr 0 \end{pmatrix},$$
where $\mbbI_E$ denotes the indicator function of the set $E.$
This yields an expansion of the residual $F(U)$ in the form
\be F(U)= \begin{pmatrix}0 \cr 2(1-p_0^2)p_0 \end{pmatrix}+\eps \begin{pmatrix}(2p_0^2-1-\eps)q_1 \cr -\beta q_1-4p_0^2p_1 \end{pmatrix}+
\eps^2\begin{pmatrix}(2p_0^2-1-\eps)q_2+4p_0p_1q_1 \cr -\beta q_2-4p_0p_1^2+L_2p_0 \end{pmatrix}+O(\eps^3). \ee
Since $p_0^2= 1$ in both domains this affords the reduction 
$$
    F(U)= \eps \begin{pmatrix}(1+\mu)q_1 \cr -\beta q_1-4p_1 \end{pmatrix}+
    \eps^2\begin{pmatrix}(1+\mu)q_2+4p_0p_1q_1 \cr -\beta q_2-4p_1^2-4p_2-2(p_1^2+q_1^2)p_0 \end{pmatrix} +O(\eps^3).
$$
The leading order outer dynamics reduces to a family of uncoupled ODEs,
$$\partial_\tau\begin{pmatrix}p_1\cr q_1\end{pmatrix} =  \begin{pmatrix}(1+\mu)q_1 \cr -\beta q_1-4p_1 \end{pmatrix}, $$
that induce exponential decay on the fast $\tau=O(1)$ timescale. Setting $u_1=0$ yields an equivalent system for $u_2$. We assume that $u_1$ and higher order outer terms are zero on the relevant time scales. Correspondingly all matching of the inner system for $i\geq 1$ is to the outer value $0.$

\subsection{Inner expansion}

The inner expansion uses the Frenet coordinates,
for which in ${\mathbb R}^n$ the scaled Laplacian takes the form
$$\eps^2\Delta=\partial_z^2+\eps \kappa_0(s)\partial_z +\eps^2 (z\kappa_1(s)\partial_z+\Delta_s)+\eps^3(\Delta_{s,1}+z^2\kappa_2(s)\partial_z) +O(\eps^4).$$
Here $\Delta_s$ is the Laplace-Beltrami operator on the interface $\Gamma$ and $\kappa_0=\sum_{j=1}^{n-1}k_j$ is the total curvature expressed in terms of the $n-1$ curvatures $\{k_1, \ldots, k_{n-1}\}$ of $\Gamma$. The higher order curvatures satisfy $$\kappa_i=(-1)^i\sum_{j=1}^{n-1}k_j^i,$$ for $i\geq 1$ see \cite{bib:DP}[eqn (2.8)] for details. In two space dimensions, with $n=2$, these relations reduce to $\kappa_i=(-1)^i\kappa_0^i.$ The inner expansion of the vector field residual 
\be \label{e:Resid} 
\tF(\tU) = \tF_0+ \eps \tF_1+ \eps^2\tF_2+\eps^3\tF_3 +O(\eps^4),\ee
requires an expansion of $\tilde U,$ 
$$
\tU = \tu_0+\eps \tu_1+\eps^2\tu_2 +O(\eps^3),
$$
where
$$\tu_0=\begin{pmatrix}\phi\cr 0\end{pmatrix}.$$
Since the higher order outer expansion is uniformly zero, the matching conditions to the outer solution devolve into requiring that each $\tu_i(s,\cdot)\in L^2(\mathbb R)$ for all values of $s$ and all $i\geq 1.$
The leading order residual has the form
\be \label{e:F0} 
\tF_0 = \begin{pmatrix}0 \cr \tE_0\phi \end{pmatrix}=0,
\ee
where we have introduced the operator
$$\tE_0:=\partial_z^3-2\phi^2+2.$$ 
This leading order residual is zero since
$\phi$ solves  \eqref{e:front-eq} which is equivalent to $\tE_0\phi=0$. 
For $i\geq 1$  the inner vector field residuals take the upper-triangular form
\be \label{e:Fi}
\tF_i= \mrL \tu_i + \tR_i,
\ee
where $\mrL$ is given in \eqref{linop}. The lower order residuals $\tR_i$ depend only upon $\tu_k$ for $k=0, \ldots i-1,$ and are given by
\be \label{e:R1}
\tR_1=  \begin{pmatrix} 0 \cr \kappa_0\partial_z\tp_0 \end{pmatrix},
\ee
\be \label{e:F2}
\tR_2=\begin{pmatrix}-\tE_1\tq_1 \cr \tE_2\tp_0+\tE_1\tp_1 \end{pmatrix},\ee
and
\be \label{e:R3}
\tR_3= \begin{pmatrix} -\tE_1\tq_2-\tE_2\tq_1
\cr \tE_1\tp_2+\tE_2\tp_1+\tE_3\tp_0 \end{pmatrix},
\ee
where
$$\tE_1=\kappa_0\partial_z-4\tp_0\tp_1,$$ $$\tE_2=z\kappa_1\partial_z+\Delta_s-2|\tu_1|^2,$$
$$\tE_3=\Delta_{s,1}+z^2\kappa_2(s)\partial_z-4\tu_1\cdot\tu_2.$$

To extract the curvature dynamics we develop a quasi-steady manifold $U$ parameterized by the interface $\Gamma$ through the scaled distance function $z$ and the curvature $\kappa_0$. These quantities evolve on the slow time $T=\eps^2 \tau$ for which $\eps^2\partial_T=\partial_\tau$. 
The chain rule gives 
\be 
\label{e:chain-rule}
D_T \tU=\frac{\partial \tU}{\partial z}\frac{\partial z}{\partial T}+\frac{\partial \tU}{\partial T}.
\ee
The normal velocity $V$ of the curve is scaled as   
$V:=-\eps^{-1} \frac{\partial z}{\partial T}.$ This affords the reduction
\be \label{e:T}
\partial_\tau \tU=\eps^2D_T \tU=-\eps V\frac{\partial \tU}{\partial z}+\eps^2\frac{\partial \tU}{\partial T}.\ee
This is further expanded in terms of the normal velocity
$$V=V_0+\eps V_1+\eps^2V_2+O(\eps^3),$$  
and the $T$ partials of $\tU$,
 $$\partial_T \tU= \eps \partial_T\tu_1+ \eps^2 \partial_T \tu_2+O(\eps^3),$$
 for which $\partial_T\tu_0=0.$
 Combining these expansions yields the inner expansion of the left-hand side of \eqref{e:FMU},
 \be \label{e:DT_exp}
  \partial_\tau\tU=-\eps V_0\partial_z \tu_0-\eps^2(V_0\partial_z\tu_1+V_1\partial_z\tu_0)-\eps^3(V_0\partial_z\tu_2+ V_1\partial_z \tu_1 +V_2\partial_z\tu_0- \partial_T\tu_1)+O(\eps^4)
\ee
Using \eqref{e:DT_exp} and \eqref{e:Fi} in \eqref{e:Resid} we match the $O(\eps)$ terms in \eqref{e:FMU}. This yields the system
\be \label{e:inner-Oeps}
-\begin{pmatrix} V_0 \cr \kappa_0 \end{pmatrix}\phi'=L \tu_1.\ee
This is an elliptic problem in $z$ for the leading order normal velocity $V_0=V_0(s)$ and $\tu_1.$ The linear operator $\mrL$ has a kernel, so Fredholm's solvability condition requires 
$$\begin{pmatrix} V_0 \cr \kappa_0 \end{pmatrix}\phi'\perp\Psi_0^\dag= \begin{pmatrix} \beta \mrD^{-1}\phi' \cr \phi' \end{pmatrix}.$$ 
This holds if the leading order normal velocity satisfies
\be V=-\alpha_0\kappa_0,
\label{e:V0}\ee
 where the curvature coefficient $\alpha_0=\alpha_0(\mu)\in\mathbb R$ satisfies
\be \label{e:alpha0}
\alpha_0:= \frac{\lVert \phi'\rVert^2}{\beta \langle \mrD^{-1}\phi',\phi'\rangle}.
\ee
This coefficient satisfies $\alpha_0>0$ so long as $\mu>0$ -- this  corresponds to motion by curvature. It is instructive to expand $\alpha_0$ in powers of $\mu$. From the relation \eqref{e:A1} we find 
\be \alpha_0=\mu\frac{\lVert \phi'\rVert^3  \lVert \psi\rVert}{\beta \langle \phi',\psi\rangle^2}+O(\mu^2).
\label{e:alpha0exp}
\ee
We deduce that $\alpha_0<0$ for $\mu$ small and negative, which implies an ill-posed motion against curvature.In the sequel this is regularized by higher order terms in the normal velocity expansion. 

\subsection{Regularization of the normal velocity}

The first step to identify higher order terms in the normal velocity is to solve the system \eqref{e:inner-Oeps} for $\tu_1$. 
The inversion formula \eqref{e:Linv-gen} applies if $\mu\neq0$, and the system can be solved directly if $\mu=0$. In either case
the correction terms have a tensor-product structure
\be 
\label{e:tu1}
\tu_1(s,z) =\kappa_0(s)\ou_1(z) =\kappa_0\begin{pmatrix}\op_1\cr \oq_1 \end{pmatrix},\ee
in terms of $s$-dependent curvature and $z$-dependent vector valued function $\ou_1$ which satisfies
\be \label{e:ou1-def}\begin{aligned}
 \ou_1 &= \begin{pmatrix}
     \mrC^{-1}\left(1-\beta \alpha_0 \mrD^{-1}\right)\phi'\cr
     \alpha_0\mrD^{-1}\phi'
 \end{pmatrix}=\alpha_0\begin{pmatrix}
     -\beta \mrC^{-1}\Pi_{\phi'}^\bot  \mrD^{-1}\phi' \cr
     \mrD^{-1}\phi'
 \end{pmatrix}, &\mu\neq 0, \\   
\ou_1&=
\begin{pmatrix} \mrC^{-1}\left(\phi'-   \frac{\Vert\phi'\Vert^2}{\langle \phi',\psi\rangle}\right)\cr
                    \frac{\Vert\phi'\Vert^2}{\beta\langle \phi',\psi\rangle} 
 \psi\end{pmatrix}=\frac{\Vert\phi'\Vert^2}{\beta\langle \phi',\psi\rangle} \begin{pmatrix}    -\beta \mrC^{-1}\Pi_{\phi'}^\bot \psi \cr  
 \psi\end{pmatrix},&\mu=0.
                        \end{aligned}
\ee
The formulas are smooth since $$\alpha_0\mrD^{-1}\phi'\to  \frac{\Vert\phi'\Vert^2}{\beta\langle \phi',\psi\rangle}\psi, $$
as $\mu\to 0.$ In particular the function $\ou_1$ is uniformly bounded as $\mu\to0$ and has even parity in $z$. 
Returning to \eqref{e:FMU}, we use \eqref{e:DT_exp} and \eqref{e:Resid} at $O(\eps^2)$. The form \eqref{e:F2} yields the balance
\be
\label{e:In-Oe2}
L\tu_2 =-\tR_2-(V_0\partial_z\tu_1+V_1\partial_z\tu_0)=\begin{pmatrix}
    \tE_1 \tq_1 -V_0\partial_z\tp_1-V_1\partial_z \tp_0 \cr -\tE_1\tp_1-\tE_2\tp_0-V_0\partial_z\tq_1
\end{pmatrix}.
\ee
The solvabilty condition for $\tu_2$ is the same as for $\tu_1,$ however all the terms on the right-hand side of \eqref{e:In-Oe2} have odd parity about $z=0$ except for $\partial_z \tp_0.$ This implies that the system is solvable for $V_1=0$. 
In two space dimensions $\kappa_1=-\kappa_0^2$  so that $\tE_2$ can be written in terms of $\kappa_0^2$. Consequently the system for $\tu_2$ has the tensor product formulation
$$ L \tu_2 = \begin{pmatrix} \omega_1\cr \omega_2 \end{pmatrix} = \kappa_0^2 \begin{pmatrix} \oomega_1 \cr \oomega_2 \end{pmatrix},$$
where we have introduced
\be\label{e:oomega1}
\oomega_1:=\partial_z(\oq_1+\alpha_0\op_1)-4\phi\op_1\oq_1,
\ee
and 
$$
\oomega_2 :=
\partial_z(\alpha_0\oq_1-\op_1)+4\phi\op_1^2+2
|\tu_1\!|^2\phi + z\phi'.
$$
The functions $\omega_1$ and $\omega_2$ have odd parity about $z=0$, in particular $\mrD^{-1}\oomega_1$ is well defined and uniformly bounded as $\mu\to 0.$
Inverting $\mrL$ we determine that
\be \label{e:tu2}
\begin{pmatrix}
    \tilde p_2\cr \tilde q_2
\end{pmatrix}= \kappa_0^2\begin{pmatrix}
    -\mrC^{-1}(\beta \mrD^{-1}\oomega_1+\oomega_2)\cr \mrD^{-1}\oomega_1
\end{pmatrix}.
\ee
This allows us to write
$$ \tu_2=\begin{pmatrix}
    \tp_2\cr \tq_2
\end{pmatrix} = \kappa_0^2(s) \begin{pmatrix}
    \op_2(z)\cr \oq_2(z)
\end{pmatrix}=\kappa_0^2\ou_2,
$$
where $\ou_2$ has odd parity in $z.$
To determine $V_2$ we proceed to the $O(\eps^3)$ matching in the inner expansion of \eqref{e:FMU}. Equating \eqref{e:R3} with  the $O(\eps^3)$ terms in \eqref{e:DT_exp}  yields
$$\begin{pmatrix}
    -V_0 \partial_z\tp_2 -V_2\partial_z\tp_0 \cr 
    -V_0\partial_z\tq_2 \end{pmatrix}
    +\partial_T\kappa_0 \ou_1 =L\tu_3+\tR_3.
 $$
The solvability conditions require that the terms without $\tu_3$ are orthogonal to $\Psi_0^\dag,$ which has even parity about $z=0$. This yields the system
\be\label{e:V2_e1}\begin{pmatrix}
    -V_0\partial_z\tp_2 -V_2 \phi' +\partial_T\kappa_0\op_1 + \tE_1 \tilde q_2+\tE_2\tilde q_1\cr
    -V_0\partial_z\tq_2 +\partial_T\kappa_0\oq_1-\tE_1\tp_2-\tE_2\tp_1-\tE_3\tp_0
\end{pmatrix}\perp \Psi_0^\dag,
\ee
to be solved for $V_2.$ From Pismen, \cite{bib:Pismen} in two space dimensions the co-moving coordinates imply the relation between normal velocity and evolution of the curvature,
\be \label{e:Pismen}
\partial_T\kappa_0=-(\Delta_s+\kappa_0^2)V=\alpha_0(\Delta_s\kappa_0+\kappa_0^3)+O(\eps^2).
\ee
This allows the left-hand side of \eqref{e:V2_e1} to be expressed as a tensor product of $\Delta_s\kappa_0$ and $\kappa_0^3$ and vector valued functions $\ov_1$ and $\ov_2$ of $z$-only dependence,
$$\left(\Delta_s \kappa_0 \ov_1 +\kappa_0^3 \ov_2 -V_2 \partial_z\ou_0 \right) \bot  \Psi_0^\dag.
$$
These $z$-only vector-valued functions take the form
\be \label{e-ov1}
\ov_1= \begin{pmatrix}
\alpha_0\op_1 +\oq_1 \cr -\op_1+\alpha_0 \oq_1
\end{pmatrix},
\ee 
and
\be \label{e-ov2}
\ov_2= \begin{pmatrix}
\alpha_0 \left(\partial_z\op_2 +\op_1\right) +\partial_z\oq_2-4\phi\op_1\oq_2-z\partial_z \oq_1-2|\ou_1|^2\oq_1 \cr
\alpha_0\left(\partial_z\oq_2+\oq_1\right) -\partial_z\op_2 +4\phi\op_1\op_2 +z\partial_z\op_1 +2|\ou_1|^2\op_1 -z\phi'+4\ou_1\cdot\ou_2\phi
\end{pmatrix}.
\ee 
Solving for $V_2$ yields the higher order corrections to the normal velocity,
\be \label{e:V2-exp}
V_2= \nu\Delta_s\kappa_0 
+\zeta \kappa_0^3.
\ee
where the coefficients are defined by
\be \label{e:nu}
\nu:= \frac{\langle \ov_1,\Psi_0^\dag\rangle}{\beta\langle \mrD^{-1}\phi',\phi'\rangle},
\ee
and
\be \label{e:zeta}
\zeta:= 
\frac{\langle \ov_2,\Psi_0^\dag\rangle}{\beta\langle \mrD^{-1}\phi',\phi'\rangle}.
\ee
Parity considerations imply that $V_3=0$ as they did for $V_1$, and hence the normal velocity has no $O(\eps^3)$ terms and
\be\label{e:V_final}
V= -\alpha_0 \kappa_0+\eps^2\left(\nu\Delta_s\kappa_0 +\zeta \kappa_0^3\right) + O(\eps^4).
\ee
Moreover the $O(\eps^4)$ terms are bounded relative to $1-\Delta_s$ and hence reflect regular perturbations.

The sign of $\nu$ is essential to the wellposedness of the normal velocity system. In particular it  requires $\nu>0.$ The definition of $\Psi_0^\dag$, \eqref{psi0mu0}, yields the formula,
\be \label{e:nu1}
\nu = \frac{ \langle \alpha_0\op_1+\oq_1,\beta \mrD^{-1}\phi'\rangle+ \langle -\op_1+\alpha_0\oq_1,\phi'\rangle}{\beta\langle \mrD^{-1}\phi',\phi'\rangle},
\ee
however $\op_1\bot\phi'$ so their inner product is zero. Using \eqref{e:ou1-def} to expand $\ou_1$ we find
\be \label{e:nu2}
\nu=-\frac{\Vert\phi'\Vert^4\beta \langle \mrC^{-1} \Pi_{\phi'}^\bot \mrD^{-1}\phi',\Pi_{\phi'}^\bot \mrD^{-1}\phi'\rangle}{\langle \mrD^{-1}\phi',\phi'\rangle^3}+ \frac{\Vert\phi'\Vert^2\Vert\mrD^{-1}\phi'\Vert^2}{\beta \langle \mrD^{-1}\phi',\phi'\rangle^2}+ \frac{\Vert\phi'\Vert^4}{\beta^3\langle \mrD^{-1}\phi',\phi'\rangle^2}.
\ee 
For $|\mu|$ small the asymptotic inverse formula
\be\label{e:Dinv} 
\mrD^{-1} \phi'= \frac{\langle \phi',\psi\rangle}{\Vert\phi'\Vert\Vert\psi\Vert}\frac{\psi}{\mu}+O(1),
\ee
shows that the second term in \eqref{e:nu2} is dominant for small $\mu$,
\be  \label{e:nu3}
\nu= \frac{\Vert\phi'\Vert^2\Vert\psi\Vert^2}{\langle \psi,\phi'\rangle^2}+O(\mu)
\ee
which is positive for $|\mu|<\mu_*,$ for $\mu_*>0$ sufficiently small, independent of $\epsilon.$ This establishes the main result \eqref{e:dsV}.

Assuming that the curvatures are uniformly bounded, the $\eps^2 \zeta\kappa_0^3$ term in \eqref{e:V_final} is asymptotically small compared to $\alpha_0 \kappa_0$ unless $\zeta$ is bounded away from zero and $\alpha_0=O(\eps^2).$ This occurs when $\mu=O(\eps^2).$ Since $\zeta=\zeta(\mu)$ is smooth in $\mu$ it remains to approximate $\zeta(0)$. The terms $\ou_1$ and $\ou_2$ are smooth in $\mu$ and in particular are bounded as $\mu\to0.$ The denominator of $\zeta$ scales like $\mu^{-1}$ as $\mu\to0$ so only the terms $\langle \ov_{21},\beta \mrD^{-1}\phi'\rangle$ can give a non-zero contribution to $\zeta(0).$ That is, for $|\mu|\ll1$ the inverse formula \eqref{e:Dinv} yields
$$ \zeta = \frac{\langle \ov_{21},\mrD^{-1} \phi' \rangle}{\langle \mrD^{-1}\phi',\phi'\rangle} +O(\mu)=\frac{\langle \ov_{21},\psi \rangle}{\langle \psi,\phi'\rangle} +O(\mu) .
$$
Since $\alpha_0\to0$ smoothly as $\mu\to0$ terms in $\ov_{21}$ containing $\alpha_0$ are also $O(\mu).$
We introduce the $z$-only reductions of the operators $\tE_1$ and $\tE_2,$
$$\oE_1:=\partial_z-4\phi\op_1,$$ $$\oE_2:=z\partial_z-2|\ou_1|^2,$$
and observe that \eqref{e:oomega1} can be written as $\oomega_1=\oE_1\oq_1+O(\mu).$ This allows $\ov_{21}$ to be expanded in the symmetric form
$$\ov_{21}=\frac{\Vert\phi'\Vert^2}{\langle\phi',\psi\rangle} \left(\oE_1\mrD^{-1}\oE_1 +\oE_2\right)\psi +O(\mu),$$
and hence
$$
\begin{aligned}
\zeta&= \frac{\Vert\phi'\Vert^2}{\langle\phi',\psi\rangle^2}
\left( \langle D^{-1} \oE_1 \psi,\oE_1^\dag \psi\rangle + \langle \oE_2\psi,\psi\rangle\right)+O(\mu),\\
&= \frac{\Vert\phi'\Vert^2}{\langle\phi',\psi\rangle^2}
\left( -\langle D^{-1}  \psi',\psi'\rangle
+ 16 \langle D^{-1}  (\phi\op_1\psi),\phi\op_1\psi\rangle
-\frac12 \Vert\psi\Vert^2 -2 \langle |\ou_1|^2,\psi^2\rangle\right)+O(\mu).
\end{aligned}
$$
The coefficient $\zeta$ at $\mu=0$, is a sum of three negative and one positive term, and hence is sign indefinite.

\section{Conclusions and Numerical Confirmation}

The length of a closed interface $\Gamma$ evolving under a normal velocity $V$ satisfies
\be\label{e:Front-length}
\partial_T |\Gamma|=\int_\Gamma V\kappa_0\,ds.
\ee
For the system \eqref{e:V_final} following an integration by parts this reduces to
\be\label{e:Fl}
\partial_T |\Gamma|=-\int_\Gamma \left(\alpha_0|\kappa_0|^2+\eps^2\nu |\nabla_s\kappa_0|^2-\eps^2\zeta|\kappa_0|^4\right)\,ds+O(\eps^4).
\ee

If the curvature is uniformly bounded by $M>0$, then the interfacial length decreases if $\alpha_0> \eps^2\zeta M^2$. In particular a circular interface $\Gamma$ with an an $O(1)$ radius $R=R_*$ is  an equilibrium if and only if $\alpha_0,\zeta>0$, $\alpha_0=O(\eps^2)$ and
\be \label{e:circ-eq}
R_*= \eps \sqrt{\frac{\zeta}{\alpha_0}}.
\ee
Conversely if $\alpha_0<0$ is $O(1)$ and the curvature is not zero, then the length of a smooth interface will grow. The curvature dynamics follow the complicated system
$$\partial_T\kappa_0 = (\Delta_s+\kappa_0^2)(\alpha_0\kappa_0 -\nu \eps^2\Delta_s \kappa_0-\eps^2\zeta \kappa_0^3)+O(\eps^4),$$
in which surface diffusion acts as a singular perturbation. Generically a simple closed interface will grow, buckle (meander), and self-intersect. Subsequent to self-intersection the front dynamics leave the regime in which the curvature driven flow was derived and the flow drives a chaotic jostling of front-type cells.
\begin{figure}[h!]
\centering
\includegraphics[width=1.9in]{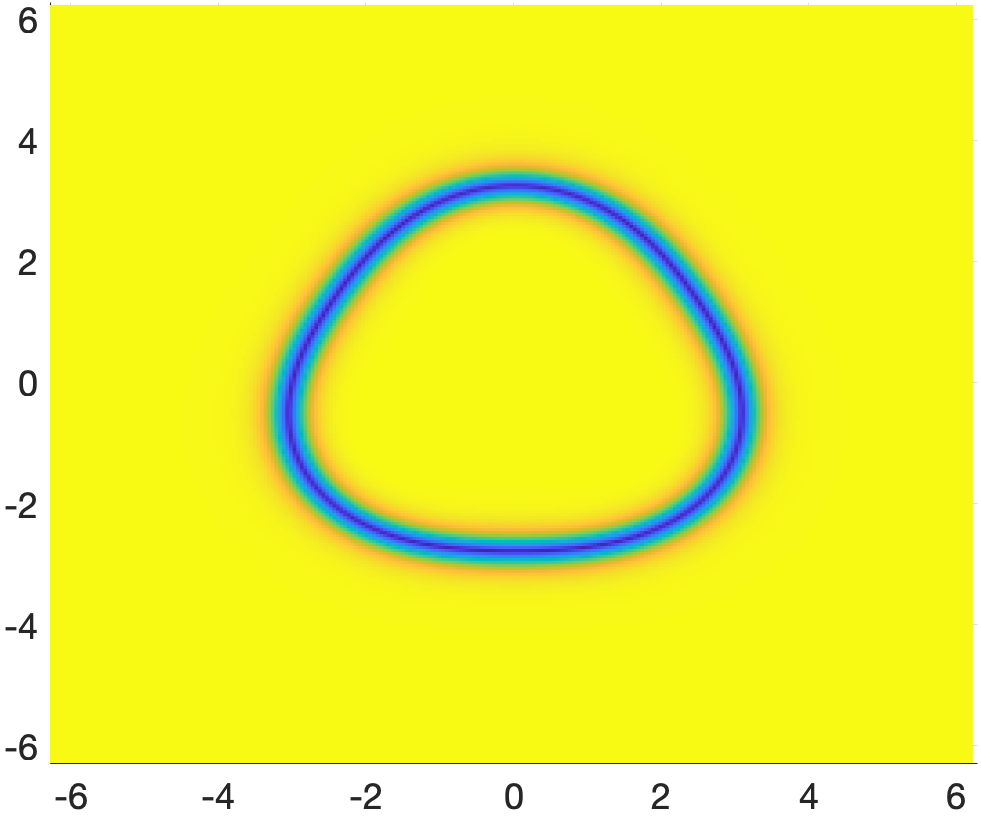}
\includegraphics[width=1.9in]{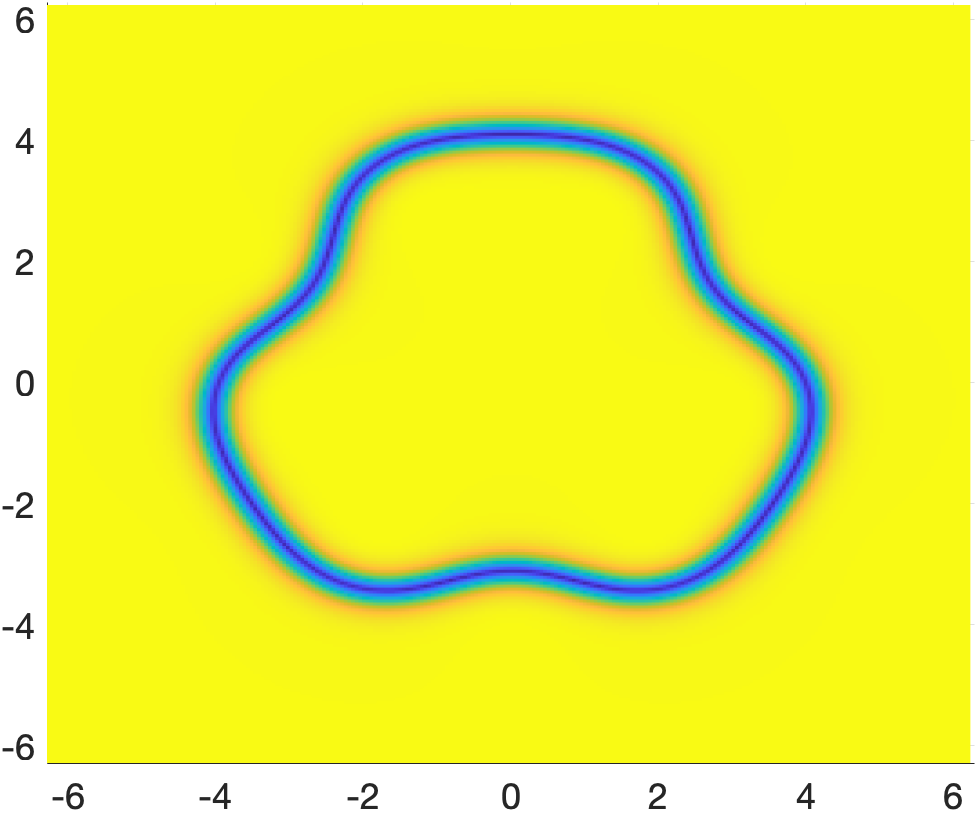}
\includegraphics[width=1.9in]{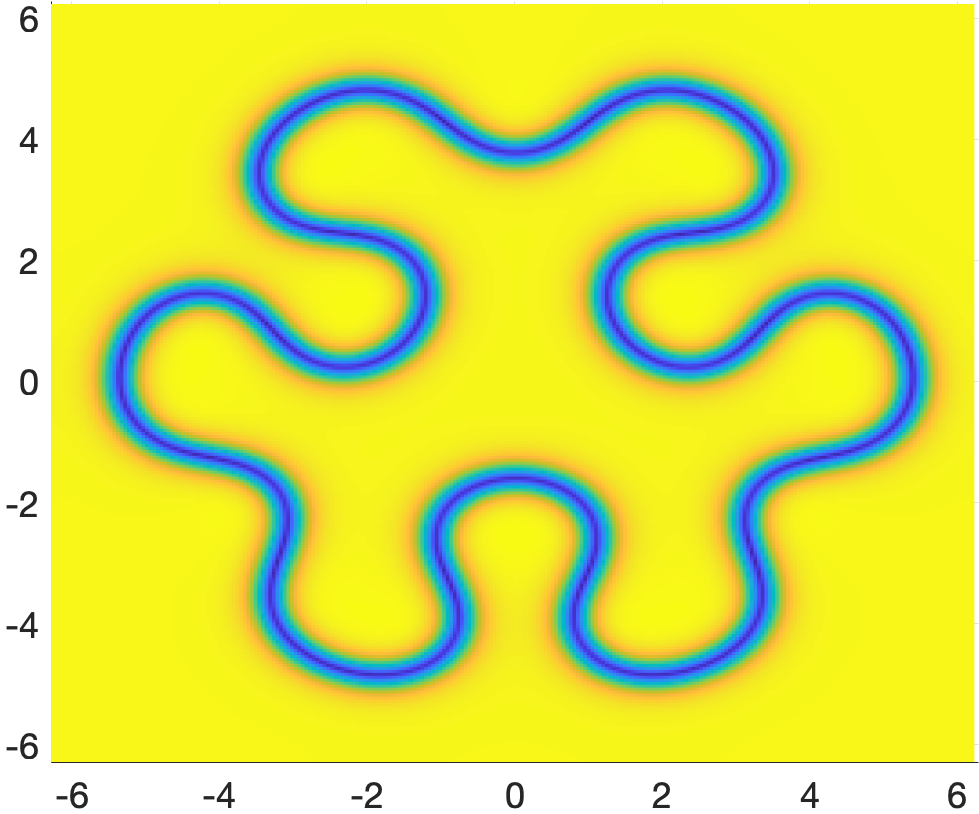}
 \put(-362, 8){\large\textbf{$\tau=0$}}
  \put(-222, 8){\large\textbf{$\tau=310$}}
   \put(-85, 8){\large\textbf{$\tau=520$}}
 \\
\includegraphics[width=1.95in]{Mu-503T0c.png}
\includegraphics[width=1.9in]{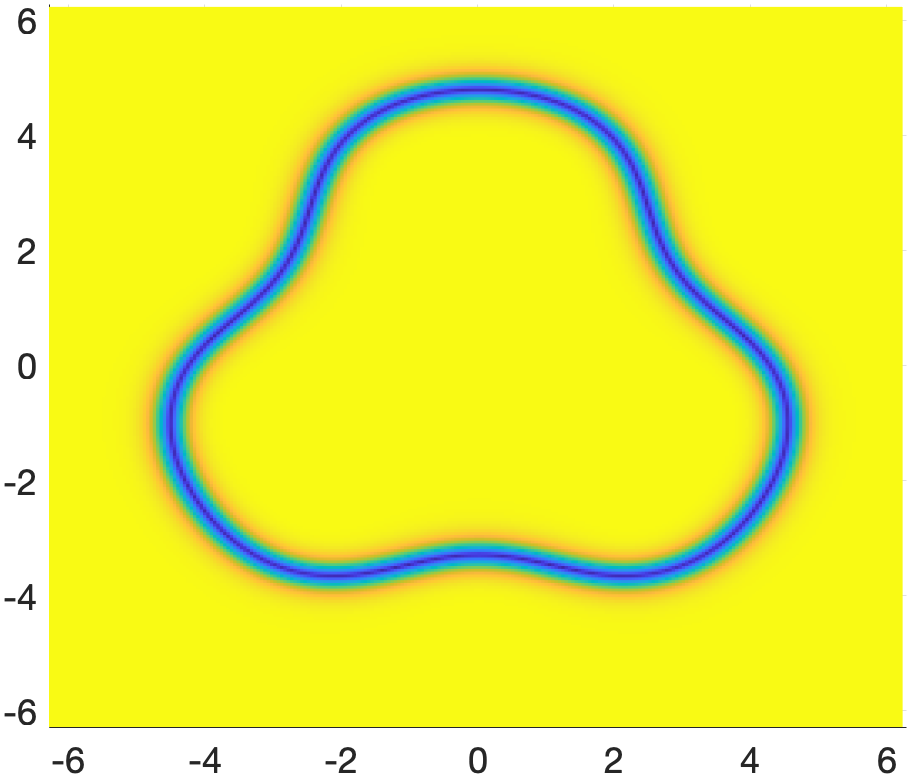} 
\includegraphics[width=1.9in]{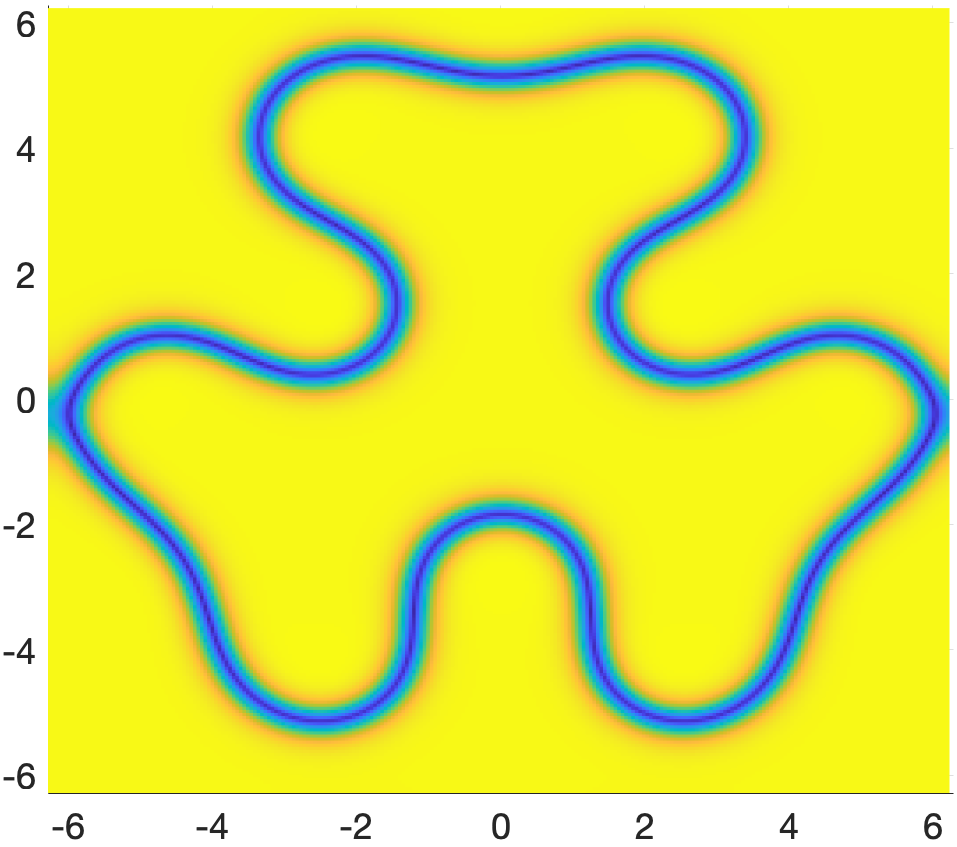} 
 \put(-362, 8){\large\textbf{$\tau=0$}}
  \put(-222, 8){\large\textbf{$\tau=1000$}}
   \put(-90, 8){\large\textbf{$\tau=1600$}}
\\
\includegraphics[width=1.95in]{Mu-503T0c.png}
\includegraphics[width=1.9in]{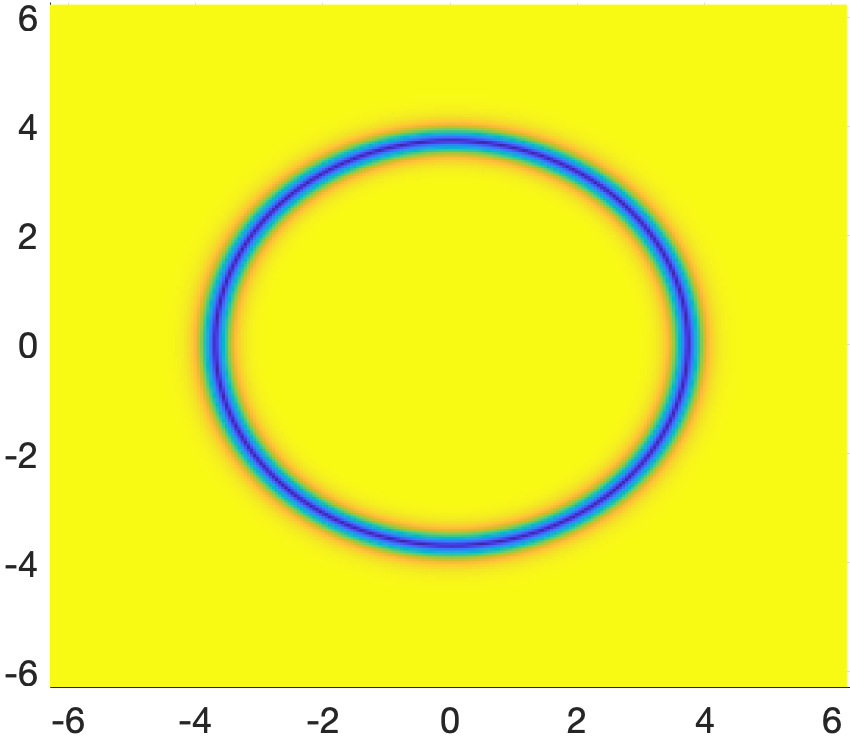} 
 \put(-222, 8){\large\textbf{$\tau=0$}}
  \put(-80, 8){\large\textbf{$\tau=10^4$}}

\caption{Plots of $|\Theta|$ as simulated by \eqref{dfpnls} for (top row) $\mu=-0.503$ and $\beta=3.50$, (middle row) $\mu=-0.208$ and $\beta=3.21$, (bottom row) $\mu=0.0405$ and $\beta=2.96.$  In all simulations $\eps=0.3$ and unscaled time $\tau$ is as indicated.}
\label{f:mu>0}
\end{figure}

These results are supported by simulations of \eqref{dfpnls} shown in Figure 1.  Images of $|\Theta|$ from the PNLS system are zero on the interface and tend to an identical constant value in both $\Omega_\pm$ domains. Each simulation starts with the same initial data,
$$\Theta_0(x)=A\tanh((|x|-r(\theta))/\eps),$$
where $A\in\mathbb C$ is the complex equilibrium of \eqref{dfpnls}, $\theta$ is the complex phase of $x$ and 
$$r(\theta)=3+\frac{1}{10}\left( \sin(3\theta) -\sin^2(7\theta)\right),$$
is a closed perturbation of a circular interface. The function $|\Theta_0|$ is depicted in the left-most image in each row of Figure 1. The dispersive ratio $\eps=0.3$ in all simulations. The top row shows the results for $\mu=-0.503$ and $\beta=3.50$ which is well into the motion against curvature regime. The interface lengthens and buckles, and self-intersects soon after the last $\tau=520$ time depicted. Subsequent evolution generates chaotic motion of front-type cells. The second row depicts the simulations for $\mu=-0.208$ and $\beta=3.21.$ This has weaker motion against curvature, the maximum curvature attained is smaller and the interface growth just yielded self-intersected (across the periodic boundary) at $\tau=1600,$ although the interface has filled the domain. The third row corresponds to $\mu=0.0405$ which is positive but smaller than $\eps^2=0.09.$ This is in the curvature driven flow regime, and the interface evolves into a circle with limiting radius $R=3.75$. The computed values $\zeta=0.576$ and $\alpha_0 = 0.0121$ yield the equilibrium radius $R_*= 3.78$, in good agreement. Computations with positive $\mu=0.159$ and $\beta=2.84$ (not shown) yield a circular interface that shrinks and approaches an $O(\eps)$ radius where it remains until $\tau=3000$ at which time the interface collapses and the function $u$ becomes spatially constant. Circular interfaces of $O(\eps)$ radius are near self-intersection and their analysis is outside the scope of this work.

\section*{Acknowledgement}
The first author acknowledges NSF support through grant DMS 2205553.


\small
\baselineskip 0pt
\begin{thebibliography}{99}

\bibitem{bib:AGJS}
  J.C. Alexander, M.G. Grillakis, C.K.R.T. Jones, B. Sandstede,
  Stability of pulses on optical fibers with phase-sensitive amplifiers,
  {\sl Z. angew. Math. Phys.} {\bf 48} (1997) 175-192.


\bibitem{bib:CP22}
Y. Chen and K. Promislow, Curve Lengthening via Regularized Motion Against Curvature from the Strong FCH Gradient Flow, {\sl J. Dynamics and Differential Equations}, (2022) DOI 10.1007/s10884-022-10178-7

  \bibitem{bib:DP}
  S. Dai and K. Promislow,
  Geometric evolution of bilayers under the functionalized Cahn–Hilliard equation{\sl Proc. R. Soc.}{\bf 469} (2013), 

 \bibitem{bib:OPO2}
  G. Iz\'us, M. Santagiustina, M. San Miguel, and P. Colet,
   Pattern formation in the presence of walk-off for a type II optical parametric oscillator, {\sl J. Opt. Soc. Am. B} {\bf 16} 1592-1596 (1999).
   
  \bibitem{bib:KP13} T. Kapitula and K. Promislow,
  {\sl Spectral and Dynamical Stability of Nonlinear Waves}, Springer, Applied Mathematical Sciences, New York, 2013.
 
  

   \bibitem{bib:Pismen}
L. M. Pismen, {\sl Patterns and interfaces in dissipative dynamics}, Springer Series in synergetics, Springer Complexity, Berlin, 2006.


 \bibitem{bib:PK}
  K. Promislow and J.N. Kutz,
  Bifurcation and asymptotic stability in the large detuning limit
  of the optical parametric oscillator, {\sl Nonlinearity} 
  {\bf 13} (2000), 675-698.
  
   \bibitem{bib:OPO4}
   C. Ropp, N. Bachelard, D. Barth, Y. Wang, and X. Zhang, Dissipative self-organization in optical space, {\sl Nature Photon} {\bf 12} 739–743 (2018).

     \bibitem{bib:OPO3}
   A. Roy, S. Jahani, C. Langrock, M. Fejer, and A. Marandi, Spectral phase transitions in optical parametric oscillators,{ \sl Nat Commun.} (2021) {\bf 12} (1) 835. 

    \bibitem{bib:OPO1}
   Majid Taki, Najib Ouarzazi, H\'elene Ward, and Pierre Glorieux, Nonlinear front propagation in optical parametric oscillators, {\sl J. Opt. Soc. Am. B} {\bf 17} 997-1003 (2000).
   
   \bibitem{bib:PRL95}
W. Zhang and J. Vinals, Secondary Instabilities and Spatiotemporal Chaos in Parametric Surface Waves,
{\sl Phys. Rev. Lett.} {\bf 74} 690  (1995).

\end{thebibliography}
\end{document}